# RESEARCH ARTICLE

# Electronic Preparatory Test for Mathematics Undergraduates: Implementation, Results and Correlations

Siri Chongchitnan,
Mathematics Institute, University of Warwick, United Kingdom.
Email: siri.chongchitnan@warwick.ac.uk

Abstract

We present a study of the implementation of the Electronic Preparatory Test for beginning undergraduates reading mathematics. The Test comprises two elements: *diagnostic* and *self-learning*. The diagnostic element identifies gaps in the background knowledge, whilst the self-learning element guides students through an upcoming material. The Test lends itself to an early identification of weak and strong students coming from a wide range of background, allowing follow-ups to be made on a topic-specific basis. The results from the Tests, collected over three years, correlate positively with end-of-year examination results. We show that such a Preparatory Test can be a better predictor of success in the first-year examination in comparison with university entry qualifications alone.

## 1. Introduction

Mathematics students starting at most universities come from a wide range of high-school qualifications and grades. In the UK, these qualifications range from A-levels to Scottish Highers and International Baccalaureates, all comprising different syllabi. Therefore, one will naturally find, in beginning undergraduates, a melting pot of abilities and backgrounds in pre-university mathematics. How, then, could university educators ensure that all students have the necessary knowledge to begin the mathematics degree?

In this paper, we will describe our work on the *Electronic Preparatory Test* (EPT) that were designed and used at the University of Hull for undergraduates reading mathematics. The EPT comprises a *diagnostic* element and a *self-learning* element. The diagnostic element ensures that students have a common baseline of background knowledge in mathematics, whilst the *self-learning* element gives a preview of the upcoming new material. Using the EPT, we were able to identify who the weaker and stronger students are, and who might need an early intervention, all within the first week of starting at university.

Pre-semester tests are being used more widely within the HE sector (see for example, Appleby, Samuels & Treasure-Jones (1997), Edwards (1997), Gillard, Levi & Wilson (2010), Sangwin (2015)). Despite many tests being implemented electronically, they often involve simple question types, or are developed on specialist electronic platforms that are not widely used and supported.

The implementation of the EPT at Hull is based on Möbius Assessment (formerly MapleTA). The Tests are online, compatible with mobile devices and graded by the Maple algebraic engine. The EPT comprises an array of interesting question types (including graph sketching and those with infinitely many correct responses). The diagnostic element is partitioned into key topics, involving a wide range of question types. The self-learning element includes an extension of school mathematics, as well as questions on university mathematics (e.g. sets and logic) that are not part of the usual school syllabi. Instead, they examine the students' ability to absorb new ideas and apply them in a self-guided way. Our method of post-Test follow-up by personal tutors will also be of interest to other HE educators.

In this work, we will begin by describing our implementation and pedagogical rationale for the EPT in detail. We will give a sample of various questions in the Test, as well as some useful electronic grading strategies. Finally, we present the results of the Tests, collected over 3 years, and study the correlation with first-year examination results.

## 2. The Electronic Preparatory Test

### 2.1 Context

The University of Hull runs a relatively small mathematics department taking in around 40 new undergraduates of mixed abilities each year. Transition from school to university mathematics was identified as a significant problem hindering progression from first to second year, and we therefore looked towards implementing a set of preparatory tests. As a result, the EPT was developed with the help of some of our own students working on creating and refining the problems during summer internships each year.

### 2.2 Question types

The EPT currently features many question types, including:

- Algebraic response;
- Free numerical response with multiple possible answers;
- Drop-down list (select the correct option);
- Tick boxes (select all correct options);
- Graph sketching.

The EPT comprises a diagnostic element and a self-learning element, which are described in detail below.

### 2.2 Diagnostic element

Students were asked to complete 5 sets of tests covering the following topics:

- **Algebra**: Binomial expansion, partial fractions, quadratic equations, simultaneous equations;
- **Numbers**: Rational numbers, integer sequences and series, inequalities;
- **Geometry**: Lines and circles in $\mathbb{R}^2$, vectors and geometric transformations;
- **Functions**: Functions and equations involving trigonometry, logarithms and the modulus;
- **Calculus**: Differentiation and integration techniques (up to integration by parts).

The diagnostic element consists of questions that probe the students' background knowledge on these topics. Here are some examples:

**Example A (from the Algebra test)**: Expand $(a-1)^4$.

**Comments**: To grade the answer, contrary to common grading strategy, it is not enough to ensure that the difference between the input and the answer evaluates to zero (otherwise, the answer $(a-1)^4$ would be marked as correct). Our grading code ensures that only $a^4 - 4a^3 + 6a^2 - 4a + 1$ or the permutation of these terms are marked as correct.

The numbers in this and other questions can be randomised, meaning that another student will not necessarily see an identical question.

**Example B (from the Calculus test):** Integrate the following function with respect to $x$. Use $C$ for the constant of integration where needed.

$$\frac{4}{x} + 1 + 3x + (2x-1)^3 + e^{5x} + \cos 2x.$$

**Comments**: All elements of the expression can be randomised. There is a small penalty for forgetting the constant of integration. The grading code can deal with answers involving permuted terms and alternative trigonometric identities.

**Example C (from the Geometry test):** Write down the equation (in the form $y = mx + c$) of a line going through $(-9,7)$ and $(5,-7)$. Sketch the line below.

**Comments:** A set of axes spanning $[-3,3]\times[-3,3]$ is given and sketching can be done by selecting two points on the line. The points can be moved with the cursor to adjust the gradient. The graphics are HTML based (there is no need to install Java or Flash) and the tolerance to the input accuracy can adjusted. Note that neither of the given points can be identified on the given axes, forcing the student to find alternative points on the lines.

## 2.3 Self-learning element

The self-learning element gives a gentle introduction to new material coming up at university. It comprises the following.

- Extension questions which are part of the five topics discussed in §2.2, but not necessarily part of school mathematics;
- A set of questions on introductory **Logic and Sets** including definitions and simple examples of sets and set operations (∩, ∪), logical implications (⟹, ⟸, ⟺), truth tables (AND, OR, NOT), logical quantifiers (∀, ∃).

Each question has self-learning material, followed by basic questions, and comments on the relevance to the upcoming modules.

Here are some sample questions.

**Example D (from the Numbers test):** A rational number is any real number that can be expressed as a fraction, $a/b$, where $a$ and $b$ are integers. A number which is not rational is called irrational.

In your first semester, you will learn a very important proof that $\sqrt{2}$ is an irrational number.

In the list below, which ones are rational numbers?

☐ $\sqrt{8}$  ☐ $-10$  ☐ 0.3333 ... (recurring)

☐ $1 + \sqrt{2}$  ☐ $\frac{3}{\sqrt{2}}$  ☐ Any rational number squared

**Comments:** Irrational numbers are not usually taught in school mathematics. Students will have to deduce the irrationality from the given fact about $\sqrt{2}$. The feedback for this question links to a YouTube video on rational numbers, and mentions the connection to the upcoming module on real analysis. The feedback also alludes to the idea of proof by contradiction, which will also be covered in an upcoming lecture.

**Example E (from the Logic and Sets test):** Logical implications are statements like:

*IF (statement A) THEN (statement B).*

For example,

*IF (Today is Tuesday) THEN (Tomorrow is Wednesday).*

*IF (I live in Hull) THEN (I live in England).*

We could use an arrow ($\Rightarrow$) to say the same thing:

*Today is Tuesday $\Rightarrow$ Tomorrow is Wednesday.*

*I live in Hull $\Rightarrow$ I live in England.*

The arrow ($\Rightarrow$) reads "implies that". You could also say it backwards as follows:

*Tomorrow is Wednesday $\Leftarrow$ Today is Tuesday.*

*I live in England $\Leftarrow$ I live in Hull.*

The arrows could go both ways if the two statements are identical (*logically equivalent*).

*Tomorrow is Wednesday $\Leftrightarrow$ Today is Tuesday.*

But of course, you can live in England without living in Hull…

Now it's your turn. In each pair of statements below, choose the symbol ($\Rightarrow$, $\Leftarrow$, $\Leftrightarrow$ from the drop-down menu) that best relates the two statements.

Let $x$ be a real number.

$x > 12$ \_\_\_\_\_\_ $x > 6$ $\qquad\qquad$ $x \leq 0$ \_\_\_\_\_\_ $x^3 \leq 0$

$x^2 = 4$ \_\_\_\_\_\_ $x = 2$ $\qquad\qquad$ $a > 0$ and $b > 0$ \_\_\_\_\_\_ $a + b > 0$

$x^2 > 4$ \_\_\_\_\_\_ $x > 2$

**Comments**: Logical implications are not part of school mathematics, but are taught within the first few weeks of the semester. It is important for students to have an early exposure to such a fundamental concept in university mathematics.

Using self-guided questions to introduce important concepts such as logical implication may seem challenging. Nevertheless, we found that the average mark for the Logic and Sets test (all of which are outside the A-level syllabus) is the highest amongst the test topics.

**Example C (from the Algebra test):** The system of equations

$$3x + 6y = 2$$
$$x + Cy = D$$

has no solution. Give values of $C$ and $D$ which ensure this is the case.

**Comments**: The grading code accepts infinitely many correct answers. The feedback for this question (regardless of whether the answers are correct) explains the interpretation of the equations as a pair of parallel lines. The feedback also explains the link to the upcoming first-year linear algebra course, in which this situation will be generalised to $\mathbb{R}^3$.

### 2.4 Rubric and follow-up

At the end of each sub-test, the students are told which answers are correct (with additional feedback, further reading, online resources and links to specific modules). The wrong answers are flagged up along with hints and explanation of common misconceptions. However, the correct answers are not shown in this case.

Students must pass all tests (the pass mark being 75%) by the end of the first week of the semester. They can repeat the tests as many times as they like (correct answers are kept and do not need to be repeated). There is no time limit for each attempt.

After the deadline, results are shared with all lecturers, and failures are flagged up to each student's personal tutor. A student who fails a particular topic will be asked to meet with their tutor to discuss the difficulties they had. They will also need to complete handwritten exercises on that topic to be submitted to the tutor by an agreed date.

## 3. Pedagogical rationale

Why do we use the EPT?

- The diagnostic element helps identify gaps in the students' background knowledge by the end of the first week of term. This means that early intervention can be administered in the form of electronic feedback and follow-up by personal tutors. Post-test follow-ups and support are essential for the success of such tests (Edwards, 1997).
- The self-learning element serves as a preview of the first-year material, hence giving students the first early exposure to university-level mathematics. It also introduces the student to the practice of pre-lecture reading used by some lecturers as part of the flipped classroom. Learning before lecture has been shown to significantly increase learning gains (Freeman et. al. (2007), Dobson (2008), Moravec et. al. (2017)).
- The EPT helps to sort students of mixed abilities by identifying who have stronger or weaker mathematics backgrounds. This helps tutors to come up with a composition of small tutorial classes that suits their teaching style (e.g. wider or narrower spread of abilities).
- In terms of confidence building, the EPT is the students' first assignment at university. Passing the EPT goes in a long way in reassuring the students that they have come to the right place, are on the right track, and are ready to start the degree. Having increased self-confidence has been associated with better examination performance (Parsons, Croft & Harrison, 2009).
- The EPT clearly sets out the rules of engagement at university. They send a clear message to the students that at university, there are strict deadlines to adhere to, as well as consequences to failing, non-completion and non-engagement.

## 4. Results

In this section, we shall describe the results of the EPT collected over 3 years (2016-2018) and their correlations with examination results.

*4.1 Correlations with examination results*

Figure 1. Data showing a positive correlation between the EPT results and average first-year exam marks (r=0.49). The data comprise results from 110 first-year students collected over 3 years (2016-2018).

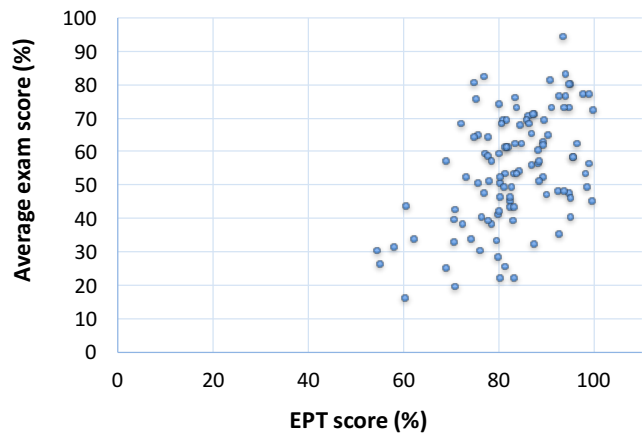

Figure 1 shows a scatter plot of the average end-of-first-year examination scores against the EPT scores for 110 students collected over 3 years (2016-2018). The examination results are calculated by averaging over all mathematics modules. Only first-attempt marks are recorded. We have excluded students who did not take all the examinations (e.g. due to medical reasons or change of course).

We calculate the correlation coefficient, $r$, defined in the usual way as

$$r_{XY} = \frac{\sum_i (x_i - \mu_X)(y_i - \mu_Y)}{\sqrt{\sum_i (x_i - \mu_X)^2 \sum_i (y_i - \mu_Y)^2}},$$

where $X$ is the average first-year examination mark and $Y$ is one of the following:

- The EPT score;
- The total entry tariff ('UCAS tariff' points based on the new system introduced in 2017, where grade A*= 56 points etc.);
- The 'maths-only' tariff, calculated using only mathematics and further mathematics grades.

The results of the calculations are shown below.

Table 1. The correlation between the first-year examination results and 3 predictors.

| **Examination VS** | $r$ |
|---|---|
| EPT score | 0.49 |
| Total entry tariff | 0.45 |
| `Maths-only' | 0.42 |

*4.2 Interpretations*

All three indicators in Table 1 correlate positively with the first-year examination results. However, the EPT appears to be the strongest correlator of first-year examination performance, followed by the total entry tariff and, rather surprisingly, the maths-only tariff.

Unlike the EPT score, the total tariff takes into account knowledge of non-mathematical subjects. Points from partial qualifications (e.g. the AS levels) and vocational training also contribute to the total tariffs of many students in our sample, thus weakening its correlation to the mathematical exam performance. Nevertheless, having a high total entry tariff could be an indicator of the ability to handle multiple responsibilities, time management and a good work ethic. Such factors are not taken into account by the maths-only tariffs, and this may explain why it is the worst of the three indicators considered.

Our results echo the findings of Lee, Harrison & Pell (2008) who found that the pre-semester diagnostics were the best predictor for first-year examination results in engineering, as well as those of Yates and James (2006) who found that entry tariffs are an often an unreliable predictor of examination results.

One interesting observation is that, if the EPT were an effective tool in helping students from a wide range of abilities succeed in the examinations, would we not expect to see a *weak* correlation between the EPT results and examination results? In other words, should all students, with the right post-test support, not achieve good exam results regardless of abilities?

To this we argue that it is unreasonable to expect the post-test support to guarantee strong examination results. The factors determining success at university are complex and wide-ranging, and certainly extend beyond background mathematical knowledge alone. Factors such as work ethic, self-confidence, motivation, conducive learning environment, peer support, health and financial stabilities, amongst a myriad of other factors, all play a role in determining success in the examination. The post-test support alone clearly cannot mitigate all these factors.

## 5. Conclusions and discussion

In this work, we have discussed the use of our Electronic Preparatory Test in mathematics as:

- an important early identification tool for weak and strong abilities, allowing for topic-specific follow-ups;
- a self-guided learn-before-lecture tool to improve learning gains;

- a predictor for first-year success, with a stronger correlation than entry tariffs;
- a first assignment to help students build confidence and understand the importance of engagement at university.

One extension idea is to implement the EPT for all returning students to ensure they retain the necessary knowledge to begin the next stage. In 2017 and 2018, we implemented an advanced version of the EPT for students who were progressing to second year, testing on key first-year ideas in real analysis, calculus, linear algebra and complex numbers. However, the sample size has so far been too small for a meaningful statistical analysis (only one-year's worth of examination results are available at the time of publication). Over the next few years, we will continue to analyse the $2^{nd}$-year EPT results and correlations with examination results and, eventually, the final degree classification.

The data for every attempt of the EPT is stored in the Möbius server. The data includes, for instance, how many times each student needed to take each sub-test before they passed and how long each attempt took. This means that it is possible to use such data to produce a better measure of the quality of each student's EPT performance (for example, by weighting the score of each test by the number of attempts needed to pass). Such an aggregate may correlate more strongly with the examination results. Another possible aggregate worth investigating is some weighted combination of the entry tariffs and the EPT results that would maximise the correlation with examination results.

Finally, we are happy to share, upon request, our EPT questions and grading code with other educators in the mathematics HE sector.

Acknowledgements: The author is grateful to Michael Hopton, Elliot Kay and Jake Vinegrad, all of whom are Hull mathematics graduates who helped with the coding and testing of the EPT.